\documentclass[10pt]{article}
\usepackage{amssymb}
\usepackage{amsmath}
\numberwithin{equation}{section}
\usepackage{graphicx}

\newtheorem{rem}{Remark}[section]
\newtheorem{thm}{Theorem}[section]
\newtheorem{lemma}[thm]{Lemma}
\newtheorem{cor}[thm]{Corollary}
\newtheorem{prop}[thm]{Proposition}

\newtheorem{example}[thm]{Example}
\newcommand\qed{\hfill\blacksquare\bigskip}

\newcommand\PP{\ensuremath{\mathbb{P}}}
\newcommand\ZZ{\ensuremath{\mathbb{Z}}}
\newcommand\CC{\ensuremath{\mathbb{C}}}

\newcommand\NN{\ensuremath{\mathbb{N}}}

\newcommand\vect[1]{\mbox{\boldmath$#1$}}
\newcommand\Abel{\widetilde{\vect{A}}}
\newcommand\abel{\vect{A}}
\newcommand\kakko[1]{\langle{#1}\rangle}
\newcommand\tee{\mathcal{T}}

\newcommand\zet[1]{\left\vert {#1} \right\vert}

\newcommand\proof{\noindent{\textbf {Proof.}}\ \,}
\begin{document}
\title{Linearisation of the $(M,K)$-reduced non-autonomous
discrete periodic KP equation}
\author{Shinsuke Iwao\\
Graduate School of Mathematical Sciences, \\
The University of Tokyo, \\
3-8-1
Komaba Meguro-ku, Tokyo \ 153-8914, Japan}
\date{}
\maketitle
\subsection*{Abstract}
The $(M,K)$-reduced 
non-autonomous discrete KP equation is linearised on the Picard group of 
an algebraic curve.
As an application,
we construct theta function solutions to the initial value problem of 
some special discrete KP equation.

\section{Introduction}\label{sec1}

The non-autonomous discrete KP equation (ndKP) is given by the formula
\cite{Willox}:
\begin{align*}
&(b(m)-c(n))\cdot f^{t+1}_{m,n}f^t_{m+1,n+1}+(c(n)-a(t))\cdot f^t_{m+1,n}f^{t+1}_{m,n+1}\\
&\ \ +(a(t)-b(m))\cdot f^{t}_{m,n+1}f^{t+1}_{m+1,n}=0, \qquad\qquad\qquad
t,m,n\in\ZZ.
\end{align*}
With constraints $a(t)=0$, $b(m)=1$, $c(n)=1+\delta_n$ and $f^t_{m,n}=f^{t-M}_{m-K,n}$,
the ndKP reduces into the following form:
$(f^{tK-mM}_n:=f^t_{m,n})$,
\[
\frac{f_{n+1}^{t+M}f_{n+1}^{t+K}}{f^{t+M+K}_{n+1}f^{t}_{n+1}}-
(1+\delta_{n+1})\frac{f^{t}_nf^{t+K}_{n+1}}{f^{t}_{n+1}f^{t+K}_n}=
-\delta_{n+1}\frac{f^{t+M+K}_nf^{t+K}_{n+1}}{f^{t+K}_nf^{t+M+K}_{n+1}}.
\]
Define $\displaystyle I_n^t:=(1+\delta_{n+1})\cdot \frac{f^t_nf^{t+K}_{n+1}}
{f^t_{n+1}f^{t+K}_n}$ and $\displaystyle
V_n^t:=\delta_{n+1}\cdot \frac{f^{t+K}_nf^{t-M+K}_{n+1}}{f^{t+K}_{n+1}f^{t-M+K}_n}$.
Then we derive the following discrete system:
$\forall n,t\in\ZZ$,
\begin{align}
&I_n^{t}=I_{n-1}^{t-M}+V_n^{t-K}-V_{n-1}^{t},\label{kdv1}\\
&V_n^{t}=\frac{I_{n}^{t-M}V_n^{t-K}}{I_n^{t}}.\label{kdv2}
\end{align}
The system (\ref{kdv1}--\ref{kdv2}) is called
\textit{$(M,K)$-reduced non-autonomous discrete KP equation} (rndKP).
The term `non-autonomous' derives from the freedom in the parameters $\delta_n$.
If we assume an extra constraint $\delta_1=\delta_2=\delta_3=\cdots$,
this system reduces to an autonomous system
(\textit{$(M,K)$-reduced autonomous discrete KP equation} (rdKP)).

In this article, we study the rndKP with the periodic boundary condition:
\begin{equation}
I_n^t= I_{n+N}^t,\quad V_n^t=V_{n+N}^t,\qquad N\in\NN.\label{kdv3}
\end{equation}
The present paper is a generalisation of the method to solve the 
generalised periodic discrete Toda equation introduced in the papers
\cite{iwao,iwao2}.
We show here that this method is also applicable to 
the quite general case of the rndKP and prove
a linearisation
theorem (theorem \ref{thm}), which illustrates the geometric information of
the discrete system.

In some special situation, theta function solutions of the initial value problem 
are constructed.
In section \ref{sec3}, we derive an explicit formula for the solutions
of the rdKP, which is a reduction of the rndKP.

{\bf Important remark}
We can assume $\mathrm{g.c.d.}(M,K)=1$ without loss of generality.
(See (\ref{kdv1}), (\ref{kdv2}).)
Aside from this, we assume $\mathrm{g.c.d.}(M+K,N)=1$ in Sections 
\ref{sec2} and \ref{sec3} by technical reason.
The general cases will be discussed in Section \ref{sec4}.
\bigskip

{\bf Notation}: For a meromorphic function $f$ over a complete curve $C$,
$(f)_0$ (resp.\,$(f)_\infty$) denotes the divisor of zeros (resp.\,poles)
of $f$. Let
$(f):=(f)_0-(f)_\infty$.
$\mathrm{Div}^d(C)$ means the set of divisors
over $C$ of degree $d$ and
$\mathrm{Pic}^d(C)$ means the quotient set defined by
$\mathrm{Pic}^d(C)=\mathrm{Div}^d(C)/(\mbox{linearly equivalent})$.
For an element $\mathcal{D}\in\mathrm{Div}^d(C)$,
$[\mathcal{D}]$ means the image of $\mathcal{D}$ under the natural map 
$\mathrm{Div}^d(C)\to\mathrm{Pic}^d(C)$.

\section{Inverse scattering method}\label{sec2}

The rndKP equation (\ref{kdv1}--\ref{kdv3}) has the following
matrix form:
\begin{equation}\label{eq2.1}
L_t(y)R_t(y)=R_{t-M}(y)L_{t-K}(y),
\end{equation}
where
\begin{gather*}
{}
L_t(y)=
\left(\begin{array}{@{\,}cccc@{\,}}
	V_1^t & 1 &  &  \\
	 & V_2^t & \ddots &  \\
	 &  & \ddots & 1 \\
	y &  &  & V_{N}^t 
\end{array}\right),\quad
R_t(y)=
\left(\begin{array}{@{\,}cccc@{\,}}
	I_1^t & 1 &  &  \\
	 & I_2^t & \ddots &  \\
	 &  & \ddots & 1 \\
	y &  &  & I_{N}^t 
\end{array}\right),
\end{gather*}
and $y$ is a complex parameter.
Let 
\begin{align}
X_t(y)&:=
L_{t-(K-1)M}(y)\cdots L_{t-2M}(y) L_{t-M}(y)L_t(y)\times\nonumber\\
&\hspace{2cm}
R_t(y)R_{t-K}(y)R_{t-2K}(y)\cdots  R_{t-(M-1)K}(y),\label{eq2.2}
\end{align}
then (\ref{eq2.1}) becomes 
\begin{equation}\label{eq2.3}
X_t(y)R_{t-MK}(y)=R_{t-MK}(y)X_{t-K}(y),
\end{equation}
or equivalently
\begin{equation}\label{eq2.3a}
L_{t-MK}(y)X_t(y)=X_{t-M}(y)L_{t-MK}(y).
\end{equation}
Because $M$ and $K$ are co-prime,
the characteristic polynomial of $X_t(y)$ does not depend on $t$.
Let 
$\widetilde{C}:=\{(x,y)\in\CC^2\,\vert\,\det{(X_t(y)-x E)}=0\}$.
Of course, $\widetilde{C}$ is also independent from $t$.
We call
the completion $C$ of $\widetilde{C}$ the \textit{spectral curve} of 
the rndKP equation.

\begin{rem}
By applying $(\ref{eq2.1})$ repeatedly,
we can transform
$(\ref{eq2.2})$ into 
\begin{align}
X_t(y)&=
R_{t-MK}(y)R_{t-(M+1)K}(y)\cdots R_{t-(2M-1)K}(y)\times\nonumber\\
&\hspace{2cm}
L_{t-(2N-1)M}(y)\cdots L_{t-(K+1)M}(y)L_{t-KM}(y).\label{eq2.2a}
\end{align}
\end{rem}

\subsection{Properties of the spectral curve}\label{sec2.1}

As a starter, we list some fundamental properties of the spectral curve $C$
in this section.
In the rest of this article, we always assume $C$ to be smooth.
Moreover, we also
assume that $\mathrm{g.c.d.}(M+K,N)=1$ in Sections \ref{sec2} and \ref{sec3}
unless otherwise is stated.

Denote the set of $N\times N$ matrices by $M_N(\CC)$ and
the subset of diagonal matrices by $\Gamma\subset M_N(\CC)$.
For a matrix $X\in M_N(\CC)$ and subsets $A,\,B\subset M_N(\CC)$,  
let 
$A+X:=\{a+X\,\vert\,a\in A\}$,
$A X:=\{a X\,\vert\,
a\in A\}$, 
$A+B:=\{a+b\,\vert\,a\in A,b\in B\}$ and 
$AB:=\{ab\,\vert\,a\in A,b\in B\}$.

Let $S$ be the $N\times N$ matrix 
$
S=:
\left(\begin{array}{@{\,}cccc@{\,}}
	0 & 1 &  &  \\
	 & 0 & \ddots &  \\
	 &  & \ddots & 1 \\
	y &  &  & 0
\end{array}\right)
$.

The polynomial $\det{(X_t(y)-x E)}$ is of degree $N$ w.r.t.{}\ $x$,
and of degree $M+K$ w.r.t.{}\ $y$.
Then the projection $p_x:C\ni (x,y)\mapsto x\in \PP$ is $(M+K):1$, and
the projection $p_y:C\ni (x,y)\mapsto y\in \PP$ is $N:1$.

Let $U_j:=(\prod_{k=1}^K{V_j^{t+k}})\cdot (\prod_{k=1}^M{I_j^{t+k}})$,
$(j\in\{1,2,\dots,N\})$.
By (\ref{kdv2}), the quantity $U_j$ is invariant under the time evolution.

\begin{prop}\label{prop2.1}
If $\mathrm{g.c.d.}(M+K,N)=1$,
the curve $C$ has the following special points:
\begin{enumerate}
\def\labelenumi{$(\mathrm{\theenumi})$}
\def\theenumi{\roman{enumi}}
\item $M$ points 
$A_j:(x,y)=\left(0,\prod_{n=1}^N{I_n^{-jK}}\right)$,\quad
$j=0,1,\dots,M-1$.
\item $K$ points 
$B_j:(x,y)=\left(0,\prod_{n=1}^N{V_n^{-jM}}\right)$,\quad
$j=0,1,\dots,K-1$.
\item
$N$ points 
$Q_j:(x,y)=
\textstyle
\left(U_j,0\right)$,\quad
$j=1,2,\dots,N$.
\item
a unique point $P:(x,y)=(\infty,\infty)$.
\end{enumerate}
\end{prop}

\proof
Let $(0,y)\in C$. 
Then we easily derive 
\begin{align*}
\prod_{j=0}^{K-1}{\det{(L_{t-jM}(y))}}\cdot
\prod_{j=0}^{M-1}{\det{(R_{t-jK}(y))}}=0,
\end{align*}
which implies
(i) and (ii).
Part (iii) follows from the fact that $L_t(0)$ and $R_t(0)$ are upper
triangular.

(iv): 
For a point $(x,y)\in C$, there exists a non-zero $N$-vector $\vect{v}(x,y)$ such that
$X_t(y)\vect{v}(x,y)=x\cdot \vect{v}(x,y)$.
Because the matrix $X_t(y)$ is contained in the subset
$(\Gamma+S)^{M+K}=\Gamma+\Gamma S+\dots+\Gamma S^{M+K-1}+S^{M+K}$,
it follows that
\begin{equation}\label{eq2.4}
(
\gamma_0+\gamma_1 S+\dots+\gamma_{M+K-1}S^{M+K-1}+ S^{M+K})\cdot\vect{v}
=x\cdot\vect{v},
\end{equation}
where $\gamma_0,\gamma_1,\dots,\gamma_{M+K-1}$ are diagonal matrices.

Define a new parameter $k$ by $y=k^{-N}$ which is assumed to be zero
near $P$. Let 
$\zeta_N$ be a $N$-th primitive root of unity.
For all $j\in\{0,1,\dots,N-1\}$,
the vector 
\[
\vect{v}_0:=((\zeta_N^jk)^{N-1}, 
(\zeta_N^jk)^{N-2},\dots,(\zeta_N^jk),\ 1)^T
\]
satisfies the formula:
$S\cdot\vect{v}_0=(\zeta_N^jk)^{-1}\cdot\vect{v}_0$.
Then,
from (\ref{eq2.4}) we obtain
$
(\zeta_N^jk)^{-M-K}\cdot\vect{v}_0=x\cdot\vect{v}_0+(\mbox{higher term})$
near $k=0$,
which implies 
$(x,y)\sim((\zeta_N^jk)^{-M-K},k^{-N})$ when $(x,y)\in C$
tends to infinity.
Because $M+K$ and $N$ are relatively prime,
we can choose an appropriate branch of $k$ around 
a unique point $P$ such that 
\[\textstyle
x=k^{-(M+K)}+\cdots,\qquad y=k^{-N}+\cdots,
\]
near $P$.
%
$\qed$

From
the proof of proposition \ref{prop2.1} one obtains
more detailed information on the point $P\in C$.
\begin{cor}\label{cor2.2}
There exists a local coordinate $k$ around $P$
such that 
\[\textstyle
x=k^{-(M+K)}+\cdots,\qquad y=k^{-N}+\cdots.
\]
$\qed$
\end{cor}
\begin{cor}\label{cor2.3}
Let $X_t(y)\vect{v}(x,y)=x\cdot \vect{v}(x,y)$.
Then, around $P$, it follows that
\[
\vect{v}(x,y)\sim(k^{N-1},k^{N-2},\dots,k,1)^T\qquad
(\mbox{up to a constant multiple}).
\]
\end{cor}
\proof Because $\mathrm{g.c.d.}(M+K,N)=1$,
the solution of the vector equation (\ref{eq2.4}) is expressed as
$\vect{v}(x,y)=(k^{N-1},k^{N-2},\dots,k,1)^T+(\mbox{higher})$
up to a constant multiple. $\qed$

\begin{rem}\label{rem2.1}
The proof of
proposition \ref{prop2.1} $(\mathrm{i})$, $(\mathrm{ii})$ implies that the set 
\[\textstyle
\left(\cup_{j=0}^{K-1}\{\prod_n{V_n^{t-jM}}\}\right)
\bigcup
\left(\cup_{j=0}^{M-1}\{\prod_n{I_n^{t-jK}}\}\right)
\]
is invariant under the time evolution.
It then follows that 
$
\{\prod{V_n^t},\prod{I_n^t}\}
=\{\prod{V_n^{t+K}},\prod{I_n^{t+M}}\}
$. 
To avoid a non-interesting solution 
$I_n^{t+M}=V_n^t$, $V_n^{t+K}=I_{n+1}^t$ of the
rndKP $(\mathrm{\ref{kdv1}}\mbox{--}
\mathrm{\ref{kdv3}})$,
we should assume the extra constraint
$
\prod_{n}{I_n^{t+M}}=\prod_{n}{I_n^{t}}\neq
\prod_{n}{V_n^{t+K}}=\prod_{n}{V_n^{t}}
$ in addition to the rndKP.
In fact, this constraint is enough
to guarantee the existence of the unique solution.
See section \ref{sec2.3}.
\end{rem}
\bigskip

Next we consider the behaviour of $Q_j$ $(j=1,2,\dots,N)$.
The position of $Q_j$ is invariant under the time evolution.
In this paper, we restrict ourselves to the following two typical cases:
\begin{center}
(a) All $Q_j$ are distinct.\qquad (b) $Q_1=Q_2=\dots=Q_N$.
\end{center}
Note that
in the case (b), the system (\ref{kdv1}--\ref{kdv3}) reduces to the rdKP.

\subsubsection*{In the case (a)}
The equation
$X_t(y)\vect{v}(x,y)=x\cdot \vect{v}(x,y)$
becomes 
\[
X_t(0)\cdot\vect{v}=U_j\cdot \vect{v},\qquad
\mbox{at } Q_j.
\]
Because $X_t(0)$ is upper triangular, the eigenvector $\vect{v}$
takes the form

\begin{equation}\label{eq2}
\vect{v}\!=\!(d_1,d_2,\dots,d_j,0,\dots,0)^T,\qquad d_j\neq 0.
\end{equation}

\subsubsection*{In the case (b)}

Let $Q:=Q_1(=Q_2=\dots=Q_N)$.
Arguments similar to those in
the proofs of corollary \ref{cor2.2}, \ref{cor2.3} 
prove the following:
\begin{prop}\label{prop2.4}
If all the points $Q_j$ coincide,
there exists a local coordinate $k$ around $Q=Q_1$ such that 
\[\textstyle
x=U_1+\cdots,\qquad y=k^N+\cdots,
\]
and the eigenvector $X_t(y)\vect{v}(x,y)=x\cdot \vect{v}(x,y)$ satisfies
\[
\vect{v}(x,y)\sim (1,k,k^2,\dots,k^{N-1})^T,\qquad 
(\mbox{up to a constant multiple}).
\]
$\qed$
\end{prop}

Let $\vect{v}(x,y)=(g_1(x,y),\dots,g_N(x,y))^T$ be an $N$-vector 
function (defined up to a constant multiples)
such that $X_t(y)\cdot \vect{v}(x,y)=x\cdot \vect{v}(x,y)$.
By the above arguments, we have:
\begin{prop}
The meromorphic function $g_j/g_{j+l}$, $(j\leq j+l\leq N)$ has:\\
$(\mathrm{i})$ $l$ zeros at $P$\\
$(\mathrm{ii})$ at least one pole at $Q_j$.
\end{prop}

Define the divisors $\mathcal{D}_1$ and $\mathcal{D}_2$
to be minimal positive divisors on $C$
such that 
\begin{gather}
(g_j/g_N)+\mathcal{D}_1\geq -(Q_j+Q_{j+1}+\dots+Q_{N-1}),\qquad \forall j, 
\label{eq2.6}\\
(g_j/g_1)+\mathcal{D}_2\geq -(j-1)P,\qquad\qquad \forall j.\label{eq2.7}
\end{gather}
These divisors were first studied in \cite{Mumford},
where it has proved that 
$\mathcal{D}_1,\,\mathcal{D}_2$ are general and
$\deg{\mathcal{D}_1}=\deg{\mathcal{D}_2}=\mathrm{genus}(C)$.

\subsection{the eigenvector mapping}\label{sec2.2}

Let $p$ be a point on a smooth curve $C$ and $k$ be a local coordinate around $p$.
For a meromorphic function $f$,
$\mathrm{ord}\,f(p)$ denotes the largest integer $r$
such that $\lim_{q\to p}\zet{k^{-r}f(q)}<+\infty$.
For a vector function
$\vect{v}(p)=(f_i)_i$, we define
$\mathrm{ord}\,\vect{v}(p):=\min_i{[\mathrm{ord}\,(f_i(p))]}$. 
\bigskip

An \textit{isolevel set} $\tee_C$ is the set of matrices $X(y)$
(eq.(\ref{eq2.2}))
associated with 
the spectral curve $C$. Let $g:=\mathrm{genus}(C)$.
Now we construct a map from $\tee_C$ to $\mathrm{Pic}^{g+N-1}(C)$ 
called the \textit{eigenvector mapping}.

Let $X=X(y)$ be an element of $\tee_C$.
If $(x,y)\in \widetilde{C}$, there exists a complex $N$-vector 
$\vect{v}(x,y)$
such that $X(y)\vect{v}(x,y)=x\,\vect{v}(x,y)$, up to a constant
multiple.
Then there exists a Zariski open subset $C^\circ$ of $\widetilde{C}$
over which the morphism $C^\circ\ni(x,y)\mapsto\vect{v}(x,y)\in\PP^{N-1}$ is 
uniquely determined.
Moreover,
for a smooth $C$, this morphism can be extended uniquely over the whole of $C$.
Denote this morphism by $\Psi_{X}:C\to\PP^{N-1}$.

The eigenvector mapping $\varphi_C:\tee_C\to\mathrm{Pic}^d(C)$\ $(d=g+N-1)$ is 
a map defined by the formula:
\begin{equation}\label{eq2.8a}
\mathcal{O}_C
(\varphi_C(X))=\Psi_{X}^\ast(\mathcal{O}_{\PP^{N-1}}(1)),
\end{equation}
where $\mathcal{O}_{\PP^{N-1}}(1)$ is the invertible sheaf
of hyperplane sections over $\PP^{N-1}$.
Note that
it is nontrivial to prove that
$\varphi_C(X)\in\mathrm{Pic}^d(C)$
(see \cite{iwao} \S 2).

Let $(X_1:X_2:\dots:X_N)$ be the homogeneous coordinate of $\PP^{N-1}$.
The eigenvector mapping illustrates the geometric 
interpretation of the general divisors $\mathcal{D}_1$ and $\mathcal{D}_2$
(section \ref{sec2.1}).
In fact, (\ref{eq2.6}) implies that
\begin{center}
$\mathcal{D}_1+Q_1+Q_2+\dots+Q_{N-1}$ is the pull-back of $\{X_N=0\}$,
\end{center}
and (\ref{eq2.7}) says
\begin{center}
$\mathcal{D}_2+(N-1)\cdot P$ is the pull-back of $\{X_1=0\}$.
\end{center}
These facts imply the following: 
$\varphi_C{(X(y))}=[\mathcal{D}_1+Q_1+Q_2+\dots+Q_{N-1}]
=[\mathcal{D}_2+(N-1)\cdot P]$.

Let $\mathfrak{d}(X(y)):=\mathcal{D}_2$.
This divisor will play an important role for constructing a tau function solution
of rdKP. See the next section.

\bigskip
\begin{rem}\label{rem2.3a}
Because $\mathcal{D}_1+Q_1+Q_2+\dots+Q_{N-1}$ and
$\mathcal{D}_2+(N-1)\cdot P$ are linearly equivalent to each other
we have 
\begin{equation}\label{eq2.9a}
(g_1/g_N)=\mathfrak{d}(X(y))+(N-1)\cdot P-\mathcal{D}_1-(Q_1+Q_2+\dots+Q_{N-1}).
\end{equation}
\end{rem}

\begin{rem}\label{rem2.3}
Let $X(y)\cdot \vect{v}(p)=x\cdot \vect{v}(p)$ and $p=(x,y)\in C$.
Equation
$(\mathrm{\ref{eq2.8a}})$ is equivalent to
$\varphi_C(X)=\left[-\sum_{p\in C}(\mathrm{ord}\,\vect{v}(p))\cdot p\right]$.
\end{rem}
\bigskip

The following theorem is essentially obtained in 
van Moerbeke, Mumford \cite{Mumford}.

\begin{thm}\label{thm2.6}
The eigenvector mapping 
$
\varphi_C:\tee_C\to\mathrm{Pic}^d(C)
$
is an embedding.
\end{thm}

%

\subsection{shift operators}\label{sec2.3}

Consider the $N\times N$ matrix $X_t(y)$ defined by (\ref{eq2.2})
and the associated spectral curve $C$.
Let $\sigma$, $\mu_K$ and $\mu_M$ be the isomorphisms on $\tee_C$ defined by:
\begin{align}
&\sigma(X_t(y)):=SX_t(y)S^{-1},\label{eq2.9}\\
&\mu_K(X_t(y)):=R_{t-(M-1)K}(y)\cdot X_t(y)\cdot \{R_{t-(M-1)K}(y)\}^{-1},\label{eq2.91}\\
&\mu_{-M}(X_t(y)):=L_{t-MK}(y)\cdot X_t(y)\cdot \{L_{t-MK}(y)\}^{-1},\label{eq2.92}
\end{align}
where
$S$ is the matrix defined in section \ref{sec2.1}.
By (\ref{eq2.3}--\ref{eq2.3a}),
we have $\mu_K(X_t)=X_{t+K}$ and $\mu_{-M}(X_t)=X_{t-M}$.
For the rndKP (\ref{kdv1}--\ref{kdv3}),
$\sigma$ is the $n$-shift operator: $n\mapsto n+1$ and 
$\mu_K$ and $\mu_{-M}$ are the $t$-shift operators: 
$t\mapsto t+K$, $t\mapsto t-M$.
Because $K$ and $M$ are co-prime,
an appropriate combination of $\mu_K$ and $\mu_M$ 
defines the unit time evolution $t\mapsto t+1$.
\bigskip

We start with the linear problem:
\begin{equation}\label{eq2.10}
X_t(y)\cdot \vect{v}(x,y)=x\cdot \vect{v}(x,y),\qquad \vect{v}(x,y)=(g_i(x,y))_{i=1}^N.
\end{equation}
This linear equation is decomposed into the following 
infinite dimensional
form:
for an infinite vector $(g_i)_{i\in\ZZ}$,
\begin{gather}
a_{i,0}\cdot g_i+a_{i,1}\cdot g_{i+1}+\dots+a_{i,M+K}\cdot g_{i+M+K}=x\cdot g_i,\qquad
(a_{i+N,j}=a_{i,j})\label{eq2.11}\\
g_{i+N}=y\cdot g_{i}\label{eq2.12}.
\end{gather}
The matrix equation (\ref{eq2.10}) can be interpreted
``(\ref{eq2.11}) with constraint (\ref{eq2.12})".
However, 
interchanging the roles of these two equations, 
\textit{i.e.},
interpreting
``(\ref{eq2.12}) with constraint (\ref{eq2.11})", we arrive at another
matrix equation:
\begin{equation}\label{eq2.13}
Y_t(x)\cdot\vect{w}=y\cdot \vect{w},\qquad\mbox{where}\quad 
\vect{w}=(g_i)_{i=1}^{M+K}. 
\end{equation}

\begin{example}
For an equation
$
\left(\begin{array}{@{\,}ccc@{\,}}
	a_1 & a_2 & 1 \\
	y & b_1 & b_2 \\
	c_2y & y & c_1
\end{array}\right)
\left(\begin{array}{@{\,}c@{\,}}
	g_1 \\
	g_2 \\
	g_3
\end{array}\right)=
x
\left(\begin{array}{@{\,}c@{\,}}
	g_1 \\
	g_2 \\
	g_3
\end{array}\right)
$,
the associated new matrix equation is:
\[
\left(\begin{array}{@{\,}cc@{\,}}
	b_2(a_1-x) & a_2b_2-b_1+x \\
	(a_1-x)(c_1-x-b_2c_2) & a_2(c_1-x)-c_2(a_2b_2-b_1-x)
\end{array}\right)
\left(\begin{array}{@{\,}c@{\,}}
	g_1 \\
	g_2
\end{array}\right)
=
y
\left(\begin{array}{@{\,}c@{\,}}
	g_1 \\
	g_2
\end{array}\right).
\]
\end{example}

We call the linear problem 
(\ref{eq2.10}) the
\textit{$x$-form} and 
the linear problem (\ref{eq2.13})
the \textit{$y$-form}.

\subsubsection{shift operators and the x-form}\label{sec2.3.1}

Due to (\ref{eq2.9}--\ref{eq2.92}),
the shift operators $\sigma$ and $\mu$ act on the eigenvector of the
$x$-form equation (\ref{eq2.10}) by:
\[\sigma:
\vect{v}\mapsto S\,\vect{v},\quad
\mu_K:
\vect{v}\mapsto \{R_{t-(M-1)K}\}\,\vect{v},\quad
\mu_{-M}:
\vect{v}\mapsto \{L_{t-MK}\}\,\vect{v}.
\]
The following lemma is easily proved:
\begin{lemma}\label{lemma2.7}
$\det{S}=(-1)^{N+1}y$,\quad 
$\det{R_{t-(M-1)K}}=\prod_n{I_n^{t-(M-1)K}}- y$,\\
$\det{L_{t-MK}}=\prod_n{V_n^{t-MK}}- y$. $\qed$
\end{lemma}

\subsubsection{shift operators and the y-form}

The $y$-form representation of the shift operators $\sigma$, $\mu_K$, $\mu_{-M}$ are
more complicated.
Let $E_1:=-(a_{1,0}-x)/a_{1,M+K}$, $E_2:=-a_{1,1}/a_{1,M+K},
\dots,E_{M+K}:=-a_{1,M+K-1}/a_{1,M+K}$.
Then (\ref{eq2.11}) becomes
$g_{M+K+1}=\sum_{i=1}^{M+K}{E_{i}g_{i}}$.

Define three new matrices $S^\ast$, $R^\ast$ and $L^\ast$ by:
\begin{gather}
S^\ast:=
\left(\begin{array}{@{\,}ccccc@{\,}}
	0 & 1 &  &  &  \\
	 & 0 & 1 &  &  \\
	 &  & \ddots & \ddots &  \\
	 &  &  & 0 & 1 \\
	E_1 & E_2 & \cdots & E_{M+K-1} & E_{M+K}
\end{array}\right),\label{eq2.15}\\
R^\ast:=
\left(\begin{array}{@{\,}ccccc@{\,}}
	I_1^- & 1 &  &  &  \\
	 & I_2^- & 1 &  &  \\
	 &  & \ddots & \ddots &  \\
	 &  &  & I_{M+K-1}^- & 1 \\
	E_1 & E_2 & \cdots & 
E_{M+K-1} & I_{M+K}^-+E_{M+K}
\end{array}\right),\label{eq2.16}\\
L^\ast:=
\left(\begin{array}{@{\,}ccccc@{\,}}
	V_1^- & 1 &  &  &  \\
	 & V_2^- & 1 &  &  \\
	 &  & \ddots & \ddots &  \\
	 &  &  & V_{M+K-1}^- & 1 \\
	E_1 & E_2 & \cdots & 
E_{M+K-1} & V_{M+K}^-+E_{M+K}
\end{array}\right),\label{eq2.16a}
\end{gather}
where $I_n^-=I_n^{t-(M-1)K}$ and $V_n^-=V_n^{t-MK}$.
The matrices
$S^\ast$, $R^\ast$ and $L^\ast$ are the $y$-form version of the matrices
$S$, $R_{t-(M-1)K}$ and $L_{t-MK}$ ``with constraint $g_{M+K+1}=\sum{E_{i}g_i}$".
The shift operators $\sigma$ and $\mu$ act on the eigenvector of the
$y$-form equation (\ref{eq2.13}) by:
\[\sigma:
\vect{w}\mapsto S^\ast\vect{w},\qquad
\mu_K:
\vect{w}\mapsto R^\ast\vect{w},\qquad
\mu_{-M}:
\vect{w}\mapsto L^\ast\vect{w}.
\]

\begin{lemma}\label{lemma2.8}
$\det S^\ast=(-1)^{M+K}\cdot (U_1-x)$,\ \ 
$\det R^\ast=\det L^\ast=(-1)^{M+K+1}\cdot x$.
\end{lemma}
\proof
The calculation is cumbersome but elementary.
We will prove this lemma in Appendix. $\qed$

\subsubsection{geometric interpretation of x-form and y-form}

Consider the projections $p_x:C\to\PP$ and $p_y:C\to\PP$ (section \ref{sec2.1}).
Recall that $p_x$ is $(M+K):1$ and $p_y$ is $N:1$.
Denote $\mathcal{F}:=\mathcal{O}_C(\varphi_C(X_t))$.

Because the $x$-form representations 
of $\sigma$, $\mu$
are independent from $x$ (section \ref{sec2.3.1}), 
for fixed $y\in\PP$ and its pre-image $p_y^{-1}(y)=\{(x_1,y),\dots,(x_{N},y)\}$,
the matrices $S$, $R_{t-(M-1)K}$ and $L_{t-MK}$
act on the vectors $\vect{v}(x_1,y),\dots,\vect{v}(x_N,y)$ simultaneously.
\footnote{If $X(y)$ has an eigenvalue $x'$ of multiplicity $m>1$,
we should choose the vectors $\vect{v}(p),\vect{v}'(p),\dots,\vect{v}^{(m-1)}(p)$,
where $\vect{v}^{(k)}(p)$ is the $k$-th differential of $\vect{v}$
with respect to the local coordinate around $p=(x',y)$.}
\footnote{Geometrically, this means that $S$, $R$ and $L$ act on the push-forward 
$(p_y)_\ast\mathcal{F}$.}

On the other hand, for generic $y$, 
the vectors $\vect{v}(x_1,y),\dots,\vect{v}(x_N,y)$ should be
linearly independent because they are eigenvectors belonging to distinct eigenvalues.
What happens if we choose $y$ such that $\det{S(y)}=0$ ?
This seemingly leads to
a contradiction, if one believes the the column vectors of the singular matrix 
$S(y)\cdot(\vect{v}(x_1,y),\dots,\vect{v}(x_N,y))$ are linearly independent.
However, realizing the fact that the eigenvectors are only
determined up to a constant,
this problem is easily solved and
we conclude that
$\sum_{i=1}^N{\mathrm{ord }\,(S(y)\vect{v}(x_i,y))}>
\sum_{i=1}^N{\mathrm{ord }\,(\vect{v}(x_i,y))}$.
\bigskip

More precisely, the statement of lemma \ref{lemma2.7} can be interpreted
as follows:
\begin{itemize}
\item 
$\sum_{i=1}^N\mathrm{ord}\,(S\vect{v}(x_i,0))=
\sum_{i=1}^N\mathrm{ord}\,(\vect{v}(x_i,0))+1$,
where $(x_i,0)\in C,\ i=1,2,\dots,N$.
\item Let $y_0:=\prod_n{I_n^{t-(M-1)K}}$.
Then 
\[\textstyle
\sum_{i=1}^N{\mathrm{ord}\,(R\vect{v}(x_i,y_0))}=
\sum_{i=1}^N{\mathrm{ord}\,(\vect{v}(x_i,y_0))} +1,\qquad
(x_i,y_0)\in C.
\]
\item Let $y_1:=\prod_n{V_n^{t-MK}}$.
Then 
\[\textstyle
\sum_{i=1}^N{\mathrm{ord}\,(L\vect{v}(x_i,y_1))}=
\sum_{i=1}^N{\mathrm{ord}\,(\vect{v}(x_i,y_1))} +1,\qquad
(x_i,y_1)\in C.
\]
\end{itemize}
\bigskip

Similar arguments in the case of the 
$y$-form representations yield the following form of
lemma \ref{lemma2.8}:
\begin{itemize}
\item 
$\textstyle
\sum_{i=1}^{M+1}{\mathrm{ord}\,(S^\ast\vect{w}(U_1,y_i))}=
\sum_{i=1}^{M+1}{\mathrm{ord}\,(\vect{w}(U_1,y_i))} +1,\qquad
(U_1,y_i)\in C
$.
\item
$
\sum_{i=1}^{M+1}{\mathrm{ord}\,(R^\ast\vect{w}(0,y_i))}=
\sum_{i=1}^{M+1}{\mathrm{ord}\,(\vect{w}(0,y_i))} +1,\qquad
(0,y_i)\in C$.
\item
$
\sum_{i=1}^{M+1}{\mathrm{ord}\,(L^\ast\vect{w}(0,y_i))}=
\sum_{i=1}^{M+1}{\mathrm{ord}\,(\vect{w}(0,y_i))} +1,\qquad
(0,y_i)\in C$.
\end{itemize}
\bigskip

Combining these data, we obtain the following proposition:
\begin{prop}\label{prop2.9}
Let $Q_1:(x,y)=(U_1,0)$ and 
\begin{align*}
A_j: &\textstyle (x,y)=(0,\prod_n{I_n^{t-(M-1)K}}),\quad
-jK\equiv t-(M-1)K \pmod{M},\\
B_i: &\textstyle (x,y)=(0,\prod_n{V_n^{t-MK}}),\quad
-iM\equiv t-MK \pmod{K}
\end{align*}
$(\mathrm{proposition\  \ref{prop2.1}})$. Then, \\[1mm]
$(\mathrm{i})$ 
$\mathrm{ord}\,(S\vect{v}(Q_1))=\mathrm{ord}\,(\vect{v}(Q_1))+1$,\quad
$(\mathrm{ii})$ 
$\mathrm{ord}\,(R\vect{v}(A_j))=\mathrm{ord}\,(\vect{v}(A_j))+1$,
$(\mathrm{iii})$ 
$\mathrm{ord}\,(L\vect{v}(B_i))=\mathrm{ord}\,(\vect{v}(B_i))+1$.
\end{prop}

\proof We prove (i).
By construction of the $x$-form and the $y$-form, we have
\[
\mathrm{ord}\,(S\vect{v}(p))=\mathrm{ord}\,(\vect{v}(p))+1
\ \Leftrightarrow\ 
\mathrm{ord}\,(S^\ast\vect{w}(p))=\mathrm{ord}\,(\vect{w}(p))+1.
\]
On the other hand,
because a regular matrix is invertible,
\[
\det{S(y)}\neq 0,\infty \ \mbox{ or } \det{S^\ast(x)}\neq 0,\infty
\quad \Rightarrow\quad 
\mathrm{ord}\,(S\vect{v}(x,y))=\mathrm{ord}\,(\vect{v}(x,y)).
\]
These facts prove proposition (i).
Clearly,
similar arguments will prove (ii) and (iii). $\qed$

\subsubsection{shift operator at the infinity point}

At $P$,
the actions $\vect{v}(P)\mapsto S\vect{v}(P)$,
$\vect{v}(P)\mapsto R\vect{v}(P)$ and $\vect{v}(P)\mapsto L\vect{v}(P)$
are directly computable.
\begin{prop}\label{prop2.10}
$(\mathrm{i})$
$\mathrm{ord}\,(S\vect{v}(P))=\mathrm{ord}\,(\vect{v}(P))-1$,\\
$(\mathrm{ii})$
$\mathrm{ord}\,(R\vect{v}(P))=\mathrm{ord}\,(\vect{v}(P))-1$,\ \
$(\mathrm{iii})$
$\mathrm{ord}\,(L\vect{v}(P))=\mathrm{ord}\,(\vect{v}(P))-1$.
\end{prop}
\proof The Proposition is readily proved by
Corollaries \ref{cor2.2} and \ref{cor2.3}. $\qed$ 

\subsection{linearisation theorem}

From the above
calculations, we obtain the linearisation theorem
representing the flow of the
rndKP equation on the Picard group of the spectral curve.

\begin{thm}\label{thm}
$(\mathrm{I})$:
Let $\mathcal{D}$ be the divisor
$
\mathcal{D}=P-Q_1.
$
Then the following 
diagram is commutative.
\[
\begin{array}{ccccc}
 &\tee_C& \to &\mbox{Pic}^d(C)& \\[2.5mm]
{\sigma}& \downarrow& &\downarrow&\hspace{-4mm}{+[\mathcal{D}]}
 \\[3mm]
  & \tee_C& \to& \mbox{Pic}^d(C)&
\end{array}.
\]
$(\mathrm{II})_1$:
Let $\mathcal{E}_{j}$ $(j=0,1,\dots,M-1)$ be the divisor
$
\mathcal{E}_j=P-A_j
$ 
and $t\equiv -(j+1)K\pmod{M}$. The following 
diagram is commutative.
\[
\begin{array}{ccccc}
 &\tee_C& \to &\mbox{Pic}^d(C)& \\[2.5mm]
{\mu_K}& \downarrow& &\downarrow&\hspace{-4mm}{+[\mathcal{E}_j]} \\[2.5mm]
  & \tee_C& \to& \mbox{Pic}^d(C)&
\end{array}.
\]
$(\mathrm{II})_2$:
Let $\mathcal{F}_{j}$ $(j=0,1,\dots,K-1)$ be the divisor
$
\mathcal{F}_j=P-B_j,
$ 
and $t\equiv -jM\pmod{K}$. The following 
diagram is commutative.
\[
\begin{array}{ccccc}
 &\tee_C& \to &\mbox{Pic}^d(C)& \\[2.5mm]
{\mu_{-M}}& \downarrow& &\downarrow&\hspace{-4mm}{+[\mathcal{F}_j]} \\[2.5mm]
  & \tee_C& \to& \mbox{Pic}^d(C)&
\end{array}.
\]
\end{thm}

\proof
The theorem follows immediately from
Remark \ref{rem2.3} and Proposition \ref{prop2.9}. $\qed$

We should note the fact that 
the position of the points $Q_j$ (proposition \ref{prop2.1}) varies under
the index shift $\sigma:n\mapsto n+1$.
To avoid confusion, we fix the rule for indexing as follows:
\textit{
Once
the points $Q_1,\dots,Q_N$ are determined, we never change their induces.
Alternatively, we define }
\begin{align*}
&\varphi_C(\sigma X(y))=\varphi_C(X(y))+[P-Q_1],\\
&\varphi_C(\sigma^2 X(y))=\varphi_C(X(y))+[P-Q_1]+[P-Q_2],\\
&\varphi_C(\sigma^3 X(y))=\varphi_C(X(y))+[P-Q_1]+[P-Q_2]+[P-Q_3],\\
&\mathit{etc.} \cdots\\
&\varphi_C(\sigma^{-1} X(y))=\varphi_C(X(y))-[P-Q_N],\\
&\varphi_C(\sigma^{-2} X(y))=\varphi_C(X(y))-[P-Q_N]-[P-Q_{N-1}],\\
&\mathit{etc.} \cdots
\end{align*}
This particular
arrangement is appropriate for our further discussion.

\begin{cor}\label{cor2.12}
Let $\mathfrak{d}(X(y))$ be the general divisor defined by
$\varphi_C(X(y))=[\mathfrak{d}(X(y))+(N-1)\cdot P]$
$(\mathrm{section}\ \mathrm{\ref{sec2.1}})$.
The divisor $\mathcal{D}_1$ in $(\mathrm{\ref{eq2.9a}})$ satisfies
$\mathcal{D}_1=\mathfrak{d}(\sigma^{-1}X(y))$.
\end{cor}
\proof
By (\ref{eq2.9a}), we obtain 
\begin{align*}
[\mathcal{D}_1]&=[\mathfrak{d}(X(y))-Q_1-\dots-Q_{N-1}+(N-1)\cdot P]\\
&=[\mathfrak{d}(\sigma^{-1}X(y))-Q_1-\dots-Q_{N-1}-Q_N+N\cdot P].
\end{align*}
By the equation $(y)=Q_1+\dots+Q_N-N\cdot P$, we conclude
$[\mathcal{D}_1]=[\mathfrak{d}(\sigma^{-1}X(y))]$. 
Because $\mathcal{D}_1$ and $\mathfrak{d}(\sigma^{-1}(X(y)))$ are
general, positive and of degree $g$,
it follows that $\mathcal{D}_1=\mathfrak{d}(\sigma^{-1}X(y))$.
$\qed$

As a conclusion of the corollary, we have
\begin{equation}\label{eq2.17}
(g_1/g_N)
=\mathfrak{d}(X)+(N-1)P-\mathfrak{d}(\sigma^{-1}X)-Q_1-\dots-Q_{N-1}.
\end{equation}
%

\section{Tau function solution of rdKP}\label{sec3}

Due to the linearisation theorem \ref{thm} and 
the injectivity of the eigenvector mapping (theorem \ref{thm2.6}),
we could say that the rndKP equation (\ref{kdv1}--\ref{kdv3}) is 
``essentially solved".
Moreover, in some fortunate case, we can construct the explicit solutions
by using the method of the \textit{Riemann theta functions}.

In the rest of the article, we assume that 
$Q_1=Q_2,\dots=Q_N(=Q)$. 
Equivalently, the rndKP reduces to the rdKP equation.
(See the paragraph after remark \ref{rem2.1}.)

Recall that we have assumed that $\mathrm{g.c.d.}(M+K,N)=1$ in the previous section.
The assumption is valid also in this section.

\subsection{construction of tau functions}

We construct a theta functional solution of rdKP equation.
As in the previous section, $X_t=X_t(y)$ denotes a square matrix defined by (\ref{eq2.2}).

Let $C$ be the (smooth) spectral curve associated with $X_t$.
Fix a symplectic basis $\alpha_1,\dots,\alpha_g;\beta_1,\dots,\beta_g$ of $C$
and the normalised holomorphic differential $\omega_1,\dots,\omega_g$ such that
$\int_{\alpha_i}{\omega_j}=\delta_{i,j}$.
The $g\times g$ matrix $\Omega:=(\int_{\beta_i}{\omega_j})_{i,j}$ is called 
the \textit{period matrix} of $C$.
For a fixed point $p_0\in C$, the \textit{Abel-Jacobi mapping} 
$\abel:\mathrm{Div}(C)\to\CC^g/(\ZZ^g+\Omega\ZZ^g)$ is a
homomorphism defined by:
\[\textstyle
\sum{Y_i}-\sum{Z_j}
\ \mapsto\  
\sum(\int_{p_0}^{Y_i}{\omega_1},\cdots,
\int_{p_0}^{Y_i}{\omega_g})-\sum(\int_{p_0}^{Z_j}{\omega_1},\cdots,
\int_{p_0}^{Z_j}\omega_g).
\]

Let us consider the universal covering $\pi:\mathfrak{U}\to C$ and
fix an inclusion $\iota:C\hookrightarrow \mathfrak{U}$.
For simplicity, we use the symbols ``$\pi$" and
``$\iota$" to express the derived maps
$\mathrm{Div}(\mathfrak{U})\to\mathrm{Div}(C)$ and
$\mathrm{Div}(C)\hookrightarrow\mathrm{Div}(\mathfrak{U})$
respectively.
Naturally, there exists a
continuous lift $\Abel:\mathrm{Div}(\mathfrak{U})\to \CC^g$
such that $\Abel\circ\iota(p_0)=0$.
For the projection $\rho:\CC^g\to\CC^g/(\ZZ^g+\Omega\ZZ^g)$,
it follows that $\rho\circ\Abel=\abel\circ\pi$.

\bigskip

Now we should define the lifted divisors 
$\mathfrak{D}(\sigma X_t),\mathfrak{D}(\mu_K X_t),\mathfrak{D}(\mu_{-M} X_t)
\in\mathrm{Div}^g(\mathfrak{U})$.
For fixed $t\in\ZZ$, assume that some lifted positive divisor $\mathfrak{D}(X_t)\in
\mathrm{Div}^g(\mathfrak{U})$ with 
$\pi(\mathfrak{D}(X_t))=\mathfrak{d}(X_t)$ is specified.
First, there uniquely exists a positive divisor 
$\mathfrak{D}(\sigma X_t)\in\mathrm{Div}^g(\mathfrak{U})$
such that:
\begin{gather}
\Abel(\mathfrak{D}(\sigma X_t))=\Abel(\mathfrak{D}(X_t)+\iota P-\iota Q),\quad
\pi(\mathfrak{D}(\sigma X_t))=\mathfrak{d}(\sigma X_t).
\end{gather}
We will consider $\mathfrak{D}(\sigma X_t)$ as the appropriately lifted divisor
of $\mathfrak{d}(\sigma X_t)$.
To choose appropriate $\mathfrak{D}(\mu_K X_t)$ and $\mathfrak{D}(\mu_{-M} X_t)$,
we have to consider the compatibility: $(\mu_K)^M+(\mu_{-M})^K=\mathrm{id}$.
On the Picard group on $C$, this reflects the equation 
\[
[(M+K)\cdot P-A_0-A_1-\dots-A_{M-1}-B_0-B_1-\dots-B_{K-1}]
=[(x)]=0.
\]
(Theorem \ref{thm} (II)).
Therefore,
we can choose $(M+K)$ points $\kappa A_0,\cdots,\kappa A_{M-1}$,
$\kappa B_0,\cdots,\kappa B_{K-1}\in\mathfrak{U}$ such that
\begin{equation}
\Abel((M+K)\cdot \iota P-\kappa A_0-\dots-\kappa A_{M-1}-
\kappa B_0-\dots-\kappa B_{K-1})=0,
\end{equation}
and $\pi(\kappa A_j)=A_j$, $\pi(\kappa B_j)=B_j$.
We now define the two divisors 
$\mathfrak{D}(\mu_K X_t)$, $\mathfrak{D}(\mu_M X_t)\in\mathrm{Div}^g(\mathfrak{U})$
by the formulas:
\begin{gather}
\Abel(\mathfrak{D}(\mu_K X_{t}))=\Abel(\mathfrak{D}(X_t)+\iota P-\kappa A_j),\quad
\pi(\mathfrak{D}(\mu_K X_{t}))=\mathfrak{d}(\mu_K X_{t}),\\
\Abel(\mathfrak{D}(\mu_{-M} X_{t}))=\Abel(\mathfrak{D}(X_t)+\iota P-\kappa B_i),\quad
\pi(\mathfrak{D}(\mu_{-M} X_{t}))=\mathfrak{d}(\mu_{-M} X_{t}),
\end{gather}
where $ t\equiv -(j+1)K\pmod{M}$, $t\equiv -iM\pmod{K}$. 
\bigskip

Let $\tau^t$ be a
holomorphic function over $\mathfrak{U}$ defined by the formula:
\begin{equation}
\textstyle
\tau^t(p)=\theta
\left(
\Abel\{\mathfrak{D}(X_t)-p-\iota\Delta\}
\right),\qquad p\in \mathfrak{U},
\end{equation}
where $\theta(\bullet)=\theta(\bullet ;\Omega)$ is the Riemann theta function
and $\Delta\in\mathrm{div}^{g-1}(C)$ is the theta characteristic divisor of $C$
(\cite{tata}, Chap.{}\,II, cor.{}\,3.11).
To avoid cumbersome notations,
we often omit the letters ``$\Abel$", ``$\iota$" and use a simpler expression
$\tau^{t}(p)=\theta(\mathfrak{D}(X_t)-p-\Delta)$,
when there is no confusion possible.

Although being defined over $\mathfrak{U}$,
$\tau^t(p)$ is considered to be a multi-valued holomorphic function over 
$C$.
By the Riemann vanishing theorem (\cite{tata}, Chap.{}\,II, thm.{}\,3.11), 
the zero divisor of
$\tau^{t}(p)$ corresponds with $\mathfrak{d}(X_t)$.
\bigskip

Let $\tau^t_+(p):=\theta(\mathfrak{D}(\sigma X_t)-p-\Delta)$,
$\tau^t_-(p):=\theta(\mathfrak{D}(\sigma^{-1} X_t)-p-\Delta)$.
Then, by theorem \ref{thm}, the function: $(\hat\sigma:=\sigma^{-1})$
\[
\Psi^t(p):=
\frac{\tau^t(p)\cdot \tau_-^{t+K}(p)}
{\tau_-^t(p)\cdot \tau^{t+K}(p)}
=\frac{\theta(\mathfrak{D}(X_t)-p-\Delta)
\cdot\theta(\mathfrak{D}(\mu_K\hat\sigma X_t)-p-\Delta)}
{\theta(\mathfrak{D}(\hat\sigma X_t)-p-\Delta)\cdot
\theta(\mathfrak{D}(\mu_K X_t)-p-\Delta)},
\]
satisfies $[(\mbox{the zeros of denominator})]=
[(\mbox{the zeros of numerator})]\in\mathrm{Pic}^{2g}(C)$ and
therefore it is a single-valued and meromorphic function over $C$.

Consider an eigenvector 
$X_t(y)
\left(\begin{array}{@{\,}c@{\,}}
	g_1^t \\
	\vdots \\
	g_N^t
\end{array}\right)
=
x
\left(\begin{array}{@{\,}c@{\,}}
	g_1^t \\
	\vdots \\
	g_N^t
\end{array}\right)
$,\ \  
$(g_i^t=g_i^t(x,y)=g_i^t(p))$.
By virtue of the equation $(g_1^t/g_N^t)
=\mathfrak{d}(X_t)+(N-1)P-\mathfrak{d}(\hat\sigma X_t)-(N-1)Q$
(\ref{eq2.9a}),
we derive the following equation from Liouville's theorem:

\begin{equation}\label{eq3.2}
\Psi^t(p)=c\times\frac{g_1^t(p)\cdot g_N^{t+K}(p)}{g_N^t(p)\cdot g_1^{t+K}(p)},
\qquad c:\mbox{constant}.
\end{equation}

Due to (\ref{eq3.2}), we can calculate some special values of $\Psi^t(p)$:
\begin{lemma}\label{lemma3.1}
If $\mathrm{g.c.d.}(M+K,N)=1$, we have
$(\mathrm{i})$
$\displaystyle
\Psi^t(P)=c$,\quad
$(\mathrm{ii})$ $\displaystyle \Psi^t(Q)=c\times \frac{I_N^{t-(M-1)K}}
{I_1^{t-(M-1)K}}$.
\end{lemma}
\proof
Because $(g_1^{t+K},\dots,g_N^{t+K})=R_{t-(M-1)K}\cdot(g_1^{t},\dots,g_N^{t})$,
we have
\begin{align*}
\Psi^t=c\times\frac{g_1^t\cdot (I_N^{t-(M-1)K}g_N^t+yg_1^t)}
{g_N^t\cdot(I_1^{t-(M-1)K}g_1^t+g_2^t)}.
\end{align*}
By corollary \ref{cor2.2}, \ref{cor2.3} and proposition \ref{prop2.4},
we easily obtain the desired result.
$\qed$

Because $\theta(\mathfrak{D}(X)-\iota Q-\Delta)=
\theta(\mathfrak{D}(X)+(\iota P-\iota Q)-\iota P-\Delta)
=\theta(\mathfrak{D}(\sigma X)-\iota P-\Delta)$,
it follows that
\[
\Psi^t(Q)=\Psi^t_+(P),\qquad\mbox{where}\quad \Psi^t_+(p)=
\frac{\tau^t_{+}(p)\cdot \tau^{t+K}(p)}{\tau^t(p)\cdot \tau_{+}^{t+K}(p)}.
\]
Then lemma \ref{lemma3.1} implies 
$I_1^{t-(M-1)K}\Psi^t_+(P)=I_N^{t-(M-1)K}\Psi^t(P)$.

Repeating the same arguments with $\Psi_+(p)$, we derive
$I_2^{t-(M-1)K}\Psi^t_{++}(P)
=I_1^{t-(M-1)K}\Psi_+^t(P)$, and inductively,
we have
\[
I_N^-\Psi^t(P)=I_1^-\Psi^t_+(P)=I_2^-\Psi^t_{++}(P)=I_3^-\Psi^t_{+++}(P)=\cdots,\quad
I_n^-=I_n^{t-(M-1)K}.
\]
Let $\Psi_n^t:=\Psi^t_{++\cdots +}(P)$ ($n$ ``$+$"s).
Finally we obtain the equations 
$\Psi_{n+N}^t=\Psi_n^t$ and
$I_n^{t-(M-1)K}\Psi_n^t=d$,
where the number $d$ does not depend on $n$.
\bigskip

Next, consider the following single-valued meromorphic function over $C$:
\[
\Phi^t(p):=
\frac{\tau_-^{t-M}(p)\cdot \tau^{t}(p)}
{\tau^{t-M}(p)\cdot \tau_-^{t}(p)}
=\frac{\theta(\mathfrak{D}(\mu_{-M}\hat\sigma X_t)-p-\Delta)
\cdot\theta(\mathfrak{D}(X_t)-p-\Delta)}
{\theta(\mathfrak{D}(\mu_{-M}X_t)-p-\Delta)\cdot
\theta(\mathfrak{D}(\hat\sigma X_t)-p-\Delta)}.
\]
Using (\ref{eq2.9a}) and Liouville's theorem,
we derive the following expression:
\begin{equation}\label{eq3.5}
\Phi^t(p)=c'\times\frac{g_N^{t-M}(p)\cdot g_1^{t}(p)}{g_1^{t-M}(p)\cdot g_{N}^{t}(p)
},\qquad c':\mbox{constant},
\end{equation}
which allows us to compute some special values of $\Phi^t(p)$.
\begin{lemma}\label{lemma3.2}
If $\mathrm{g.c.d.}(M+K,N)=1$, we have
$(\mathrm{i})$ $\displaystyle
\Phi^t(P)=c'$,\quad
$(\mathrm{ii})$ $\displaystyle
\Phi^t(Q)=c'\times\frac{V_{N}^{t-KM}}{V_1^{t-KM}}$.
\end{lemma}
\proof By $(g_1^{t-M},\dots,g_N^{t-M})=L_{t-KM}\cdot
(g_1^t,\dots,g_N^t)$,
it follows that
\[
\Phi^t=c\times \frac{(V_N^{t-KM}g_1^{t}+yg_1^{t})\cdot g_1^{t}}
{(V_1^{t-KM}g_N^{t}+g_2^{t})\cdot g_N^{t}}.
\]
By virtue of
corollaries \ref{cor2.2}, \ref{cor2.3} and proposition \ref{prop2.4},
we obtain the desired result.
$\qed$

Due to $\Phi^t(Q)=\Phi^t_+(P)$ and lemma \ref{lemma3.2},
we have $V_1^{-}\Phi_+^t(P)=V_{N}^{-}\Phi^t(P)$, which implies
\[
V_{N}^-\Phi^t(P)\!=\!V_1^-\Phi_+^t(P)
\!=\!V_2^-\Phi_{++}^t(P)\!=\!V_3^-\Phi_{+++}^t(P)=\cdots,\quad V_n^-\!=\!V_n^{t-KM}.
\]
Let $\Phi_{n}^t:=\Phi_{++\cdots+}^t(P)$ ($n$ ``$+$"s).
Therefore we obtain 
$\Phi_{n+N}^t=\Phi_n^t$ and
$V_n^{-}\Phi_n^t=d'$,
where the number $d'$ does not depend on $n$.
\bigskip

Define 
$\tau_{-1}^t:=\tau_-^t(\iota P)$, 
$\tau_0^t:=\tau^t(\iota P)$, 
$\tau_1^t:=\tau_{+}^t(\iota P),
\cdots,
\tau_{n}^t:=\tau^t_{++\dots+}(\iota P)$ ($n$ ``$+$"s).
By the above arguments,
$I_n^t$ and $V_n^t$ have following expressions:

\begin{equation}\label{eq3.3}
I_n^t=d\times
\frac{\tau_{n-1}^{t+(M-1)K}\cdot \tau_n^{t+MK}}
{\tau_n^{t+(M-1)K}\cdot \tau_{n-1}^{t+MK}},\qquad
V_n^t=d'\times\frac{\tau_{n}^{t+(K-1)M}\cdot \tau_{n-1}^{t+KM}}
{\tau_{n-1}^{t+(K-1)M}\cdot\tau_n^{t+KM}}.
\end{equation}

\subsection{solution of rdKP}

For $g$-dimensional vectors $\vect{a}$ and $\vect{b}$, 
$\langle \vect{a},\vect{b} \rangle$ denotes
$\vect{a}^T \vect{b}\in\CC$.

Due to the periodicity $\mathfrak{d}(\sigma^N X_t)=\mathfrak{d}(X_t)$,
there exist integer vectors $\vect{n},\,\vect{m}\in\ZZ^g$ such that 
$
\Abel(N(\iota P-\iota Q))=\vect{n}+\Omega\vect{m}
$.
Considering the definition of the Riemann theta function 
(see \cite{tata}, \S II.1, for example),
we have 
\[
\tau^t_{n+N}=\tau_{n}^t\times\exp(-2\sqrt{-1}\,\pi\cdot
\langle\vect{m},\vect{z}\rangle-\sqrt{-1}\,\pi\cdot
\langle \vect{m},\Omega\vect{m}\rangle),
\]
where $\vect{z}=\Abel(\mathfrak{D}(\sigma^{n}X_t)-\iota P-\Delta)$.
By (\ref{eq3.3}), we have
\begin{align}
I_1^tI_2^t\cdots I_N^t&=d^N\times\frac{\tau_0^{t+(M-1)K}\cdot\tau_{N}^{t+MK}}
{\tau_{N}^{t+(M-1)K}\cdot\tau_0^{t+MK}}\nonumber\\
&=d^N\times \exp(-2\sqrt{-1}\,\pi\cdot
\langle  \vect{m}, \Abel(\iota P-\kappa A_j) \rangle),\label{eq3.6}\\
V_1^tV_2^t\cdots V_N^t&={d'}^N\times
\frac{\tau_{N}^{t+(K-1)M}\cdot\tau_{0}^{t+KM}}
{\tau_0^{t+(K-1)M}\cdot\tau_{N}^{t+KM}}\nonumber\\
&\hspace{-1.5cm}
={d'}^N\times\exp(-2\sqrt{-1}\,\pi\cdot\langle
\vect{m},\Abel(\iota P-\kappa B_i)
\rangle),
\label{eq3.7}
\end{align}
where $t\equiv -jK\pmod{M}$ and $t\equiv -iM\pmod{N}$.
Recall $\prod_n{I_n^{t+M}}=\prod_n{I_n^t}$ and $\prod_n{V_n^{t+K}}=\prod_n{V_n^t}$,
which imply that 
$d$ depends on $t\pmod{M}$ and that $d'$ depends on $t\pmod{K}$.
Finally we obtain the conclusion:

\begin{thm}
On condition that $\mathrm{g.c.d.}(M+K,N)=1$,
$(\ref{eq3.3}\mbox{--}\ref{eq3.7})$ solves
the rdKP equation $(\ref{kdv1}\mbox{--}\ref{kdv3})$ with
constraint 
\[
U_1=U_2=\dots=U_N.
\]
\end{thm}


\section{The general cases}\label{sec4}

In the previous sections, we have assumed that $M+K$ and $N$ are relatively 
prime.
This is the time to delete the assumption and discuss the general cases.

Unfortunately, the method which we have established in this paper cannot be applied
to the general cases.
For example, when $M=K=1$, $N=2$, the defining polynomial
of the spectral curve (Section \ref{sec2}) is
\[
\det{(X_t(y)-xE)}=y^2-y(2x+U_1)+x^2-U_2x+U_3,
\]
where $U_1=I_1^tI_2^t+V_1^tV_2^t$, 
$U_2=V_1^tI_1^t+V_2^tI_2^t$, $U_3=I_1^tI_2^tV_1^tV_2^t$.
Of course, these $U_1$, $U_2$ and $U_3$ are conserved quantities of the 
discrete reduced KP system (\ref{kdv1}), (\ref{kdv2}).
However, there exists another hidden independent conserved quantity of the system.
In fact, $I_1^t+I_2^t+V_1^t+V_2^t$ is invariant under the time evolution $t\mapsto t+1$
and is independent from $U_1$, $U_2$ and $U_3$.
This means that the spectral curve fails to reflect faithfully the 
information of the system.

Therefore, we should construct a new method for the general cases.
Now we prove that every reduced KP equations can be traced to the case 
$\mathrm{g.c.d.}(M+K,N)=1$.
Denote by $\mathrm{KP}_{M,K,N}$ the reduced discrete KP equation (\ref{kdv1}--\ref{kdv2})
associated with the positive integers $M$, $K$ and $N$.

Let 
$\left\{
\begin{array}{l}
\Lambda:=\ZZ_{\geq 0}\setminus\bigcup_{k=0}^\infty{\{kM,kM+1,\dots kM+K-1\}}\\
\Xi:=\ZZ_{\geq 0}\setminus\bigcup_{k=0}^\infty{\{kK,kK+1,\dots kK+M-1\}}
\end{array}
\right.
$.
Note that 
$K<M\ \Leftrightarrow \Lambda\neq\emptyset$ and
$M<K\ \Leftrightarrow \Xi\neq\emptyset$.

\begin{prop}
$(\mathrm{i})$
Suppose $K<M$.
Define the initial values $I_n^0,\dots,I_n^{K-1}:=\zeta+o(\zeta)$,
$(\zeta\to\infty,\forall n)$ for some complex parameter $\zeta$.
If $\{I_n^t,V_n^t\}_{n,t\in\ZZ}$ is a solution of $\mathrm{KP}_{M,K,N}$,
then the sequence $\{I_n^t,V_n^t\}_{n\in\ZZ,t\in\Lambda}$ converges to a solution of
$\mathrm{KP}_{M-K,K,N}$ when $\zeta\to\infty$.
\end{prop}

\proof Order the elements of $\Lambda$ as $\Lambda=\{t_1<t_2<t_3<\cdots\}$.
Note that $t_{s-M+K}=t_s-M$ and 
\[
t_s-t_{l}\equiv 0 \pmod{K}\ \Leftrightarrow\ l\equiv 0\pmod{K}, \qquad
(\because\ t_{s}-t_{s-1}=1,\mbox{ or }K+1).
\]
Especially, we have $t_s-t_{s-K}=kK\,\Rightarrow t_s-K,t_s-2K,\dots,t_s-(k-1)K\not\in\Lambda$.

To prove the statement, it is sufficient to say 
\[
\left\{
\begin{array}{l}
I_n^{t_s}=I_{n-1}^{t_{s-M+K}}+V_n^{t_{s-K}}-V_{n-1}^{t}+o(1)\\
\displaystyle V_n^{t_s}=\frac{I_{n}^{t_{s-M+K}}V_n^{t_{s-K}}}{I_n^{t}}\cdot (1+o(1))
\end{array}
\right.,
\]
or equivalently, 
\begin{equation}\label{eq4.1}
\left\{
\begin{array}{l}
I_n^{t}=I_{n-1}^{t-M}+V_n^{t-kK}-V_{n-1}^{t}+o(1)\\
\displaystyle V_n^{t}=\frac{I_{n}^{t-M}V_n^{t-kK}}{I_n^{t}}\cdot (1+o(1))
\end{array}
\right.,\ \ 
(k\in\ZZ_{>0}; t,t-kK\in\Lambda).
\end{equation}

By (\ref{kdv1}--\ref{kdv2}) and remark \ref{rem2.1}, we have
\[
I_n^{t-M}=\zeta+o(\zeta),\ (\forall n)\ \Rightarrow\ 
\left\{
\begin{array}{l}
I_n^t=\zeta+o(\zeta),\ (\forall n)\\
V_n^t=V_n^{t-K}+o(1),\ (\forall n)
\end{array}
\right..
\]
In our situation, it follows that $t\not\in\Lambda\,\Rightarrow\,
V_n^t=V_n^{t-K}+o(1),\ (\forall n)$.
Using (\ref{kdv1}--\ref{kdv2}) again, we can conclude (\ref{eq4.1}) soon.
$\qed$

\addtocounter{thm}{-1}

Similarly, we have:
\begin{prop}
$(\mathrm{ii})$
Suppose $M<K$.
Define the initial values $V_n^0,\dots,V_n^{M-1}$ $:=\zeta+o(\zeta)$,
$(\zeta\to\infty,\forall n)$ for some complex parameter $\zeta$.
If $\{I_n^t,V_n^t\}_{n,t\in\ZZ}$ is a solution of $\mathrm{KP}_{M,K,N}$,
then the sequence $\{I_n^t,V_n^t\}_{n\in\ZZ,t\in\Xi}$ converges to a solution of
$\mathrm{KP}_{M,K-M,N}$ when $\zeta\to\infty$. $\qed$
\end{prop}

\begin{example}
The reduced discrete KP equation with $M=K=1$, $N=2$ can be traced to $M=2$, $K=1$, $N=2$.

Let 
$
L_1=
\left(\begin{array}{@{\,}cc@{\,}}
	V_1^1 & 1 \\
	y & V_2^1
\end{array}\right)
$,
$
R_0=
\left(\begin{array}{@{\,}cc@{\,}}
	\zeta & 1 \\
	y & \zeta
\end{array}\right)
$,
$
R_1=
\left(\begin{array}{@{\,}cc@{\,}}
	I_1^1 & 1 \\
	y & I_2^1
\end{array}\right)
$, and $X_1:=L_1R_1R_0$.
The defining function of the spectral curve is 
\begin{align*}
\det{(X_1(y)-xE)}&=-y^3+y^2(\zeta^2+U_1)
-y\{x(2\zeta+U_4)+\zeta^2U_1+U_3\}\\
&\hspace{1cm}+x^2-\zeta U_2x+\zeta^2 U_3,
\end{align*}
where $U_4=I_1^1+I_2^1+V_1^1+V_2^1$. Note that $U_4$ is the hidden conserved quantity
of $\mathrm{KP}_{1,1,2}$.
If $\{I_n^t,V_n^t\}_{n,t}$ is a solution of $\mathrm{KP}_{2,1,2}$, the sequence
\[
\lim_{\zeta\to\infty}{I_n^1},\lim_{\zeta\to\infty}{I_n^3},\lim_{\zeta\to\infty}{I_n^5},\dots;
\lim_{\zeta\to\infty}{V_n^1},\lim_{\zeta\to\infty}{V_n^3},\lim_{\zeta\to\infty}{V_n^5},\dots
\]
is a solution of $\mathrm{KP}_{1,1,2}$.
\end{example}

\subsection*{Acknowledgement}

The author is very grateful to Professor Tetsuji Tokihiro and Professor Ralph
Willox for helpful comments on this paper.
This work was supported by
KAKENHI 09J07090.

\appendix
\section{Proof of lemma \ref{lemma2.8}}

Let $S^\ast$, $R^\ast$ and $L^\ast$ 
be the matrices defined by (\ref{eq2.15}--\ref{eq2.16a}).
We first calculate the coefficients $a_{i,j}$ of the equation 
(\ref{eq2.11}--\ref{eq2.12}):
\[
a_{i,0}\cdot g_i+a_{i,1}\cdot g_{i+1}+\dots+a_{i,M+K}\cdot g_{i+M+K}=x\cdot g_i,\quad
g_{i+N}=yg_i,
\]
which is equivalent to the formula $X(y)\cdot(g_1,\dots,g_N)=x\cdot (g_1,\dots,g_N)$.

Denote the vector $(g_1,\dots,g_N)$ by $(g_i)_i$ simply.
By equation (\ref{eq2.2}): 
\begin{align*}
X_t(y)&:=
L_{t-(K-1)M}(y)\cdots L_{t-2M}(y) L_{t-M}(y)L_t(y)\times\\
&\hspace{2cm}
R_t(y)R_{t-K}(y)R_{t-2K}(y)\cdots  R_{t-(M-1)K}(y),
\end{align*}
we have
\begin{align*}
&X_t\cdot(g_i)_i=L_{t-(K-1)M}\cdots L_{t-M}L_t
R_tR_{t-K}\cdots R_{t-(M-1)K}\cdot (g_i)_i\\
&=L_{t-(K-1)M}\cdots L_{t-M}L_t
R_tR_{t-K}\cdots R_{t-(M-2)K} \cdot(I_i^{t-(M-1)K}g_i+g_{i+1})_i\\
&=L_{t-(K-1)M}\cdots L_tR_t\cdots R_{t-(M-3)K}\cdot 
(I_i^{t-(M-2)K}\{I_i^{t-(M-1)K}g_i+g_{i+1}\}\\
&\hspace{7cm} +\{I_{i+1}^{t-(M-1)K}g_{i+1}+g_{i+2}\})_i\\
&=\cdots.
\end{align*}

Let 
$\mathcal{X}_l$ be the set of sequences of letters $s$ and $m$
of length $l$.
For example, $\mathcal{X}_0=\emptyset$,
$\mathcal{X}_1=\{s,m\}$, $\mathcal{X}_2=\{ss,sm,ms,mm\},\cdots$.
Denote $\mathcal{X}:=\cup_l\mathcal{X}_l$.
Consider a map $\kakko{\cdot} :\mathcal{X}\to \CC$ defined by: 
$\kakko{s}:=1$, $\kakko{m}:=I_i^{t-(M-1)K}$ and
\begin{gather*}
\kakko{s\chi}:=\sigma\kakko{\chi},\quad\kakko{m\chi}:=I_i^{t-(M-l)K}\cdot
\kakko{\chi},\qquad 
\chi\in\mathcal{X}_{l-1}\ (1< l\leq M),\\
\kakko{s\chi}:=\sigma\kakko{\chi},\quad\kakko{m\chi}:=
V_i^{t-(l-1)M}\cdot\kakko{\chi},\qquad 
\chi\in\mathcal{X}_{M+l-1}\ (1\leq l< K),\\
\textstyle
\mbox{where}\qquad\quad \sigma:\prod{I_j^s}\cdot\prod{V_j^s}\mapsto 
\prod{I_{j+1}^s}\cdot\prod{V_{j+1}^s}
\mbox{ is the index shift operator.}
\end{gather*}
This notation allows us to express the vector $X_t\cdot (g_i)_i$.
\begin{lemma}\label{lemmaa.1}
Let 
\[
\mathcal{X}_{l,k}:=
\{\chi\in\mathcal{X}_{l}\,\vert\,\mbox{The number of `$s$' contained in $\chi$ is $k$}\},
\]
and
$a_{i,k}:=\sum_{\chi\in\mathcal{X}_{M+K,k}}{\kakko{\chi}}$.
Then 
\[
X_t\cdot (g_i)_i=(a_{i,0}\cdot g_i+a_{i,1}\cdot g_{i+1}+\dots+a_{i,M+K}
\cdot g_{i+M+K})_i.
\]
\end{lemma}
\proof
Let $\Xi_l:=R_{t-(M-l)K}(y)$ $(l=1,2,\dots,M)$,
$\Xi_{M+l}:=L_{t-(l-1)M}$ $(l=1,2,\dots,K)$ and 
$\xi_l:=I_i^{t-(M-l)K}$ $(l=1,2,\dots,M)$,
$\xi_{M+l}:=V_i^{t-(l-1)M}$ $(l=1,2,\dots,K)$.
Assume 
\[\textstyle
\Xi_l\Xi_{l-1}\cdots \Xi_1\cdot (g_i)_i=
(\sum_k\sum_{x\in \mathcal{X}_{l,k}}{\kakko{x}}\cdot g_{i+k})_i.
\]
Then we have
\begin{align*}
\textstyle
&\Xi_{l+1}\Xi_{l}\cdots \Xi_1\cdot (g_i)_i\\
&\textstyle
=\Xi_{l+1}\cdot (\sum_k\sum_{\chi\in \mathcal{X}_{l,k}}{\kakko{\chi}}\cdot g_{i+k})_i\\
&\textstyle
=(\xi_{l+1}\cdot \{\sum_k\sum_{\chi\in \mathcal{X}_{l,k}}{\kakko{\chi}}\cdot g_{i+k}\}
+\sigma\{\sum_k\sum_{\chi\in \mathcal{X}_{l,k}}{\kakko{\chi}}\cdot g_{i+k}\}
)_i\\
&\textstyle
=(\sum_k\sum_{\chi\in \mathcal{X}_{l,k}}{\kakko{m\chi}}\cdot g_{i+k}+
\sum_k\sum_{\chi\in \mathcal{X}_{l,k}}{\kakko{s\chi}}\cdot g_{i+k+1})_i\\
&\textstyle
=(\sum_k\sum_{\chi'\in \mathcal{X}_{l+1,k}}{\kakko{\chi'}}\cdot g_{i+k})_i.
\end{align*}
By induction, we obtain 
$X_t\cdot (g_i)_i=
\Xi_{M+K}\Xi_{M+K-1}\cdots \Xi_1\cdot (g_i)_i=
(a_{i,0}\cdot g_i+a_{i,1}\cdot g_{i+1}+\dots+a_{i,M+K}
\cdot g_{i+M+K})_i$. $\qed$

In particular, we have $a_{i,M+K}=\kakko{ss\cdots s}=1$.

The following lemma is given for use in the calculations further below.
\begin{lemma}\label{lemmaa.2}
Let $\chi\in\mathcal{X}_{l,k}$.
Then $\kakko{\chi m}=\kakko{\chi s}\cdot I_{i+k}^{t-(M-1)K}$.
\end{lemma}
\proof
Let $I_i^-:=I_i^{t-(M-1)K}$ and
assume $\kakko{\chi'm}=\kakko{\chi's}\cdot I_{i+h}^{-}$ for
$\chi'\in\mathcal{X}_{l-1,h}$\, $(l-1\geq h)$.
If $\chi=s\chi'$ $(y\in\mathcal{X}_{l-1,k-1})$, we have 
$\kakko{\chi m}=\kakko{s\chi'm}=\sigma(\kakko{\chi'm})
=\sigma(\kakko{\chi's}\cdot I_{i+k-1}^{-})
=\kakko{s\chi' s}\cdot I_{i+k}^{-}=\kakko{\chi s}\cdot I_{i+k}^{-}$.
If $\chi=m\chi'$ $(\chi'\in\mathcal{X}_{l-1,k})$, $\kakko{\chi m}=\kakko{m\chi'm}
=\xi_{l}\cdot\kakko{\chi'm}=\xi_l\cdot\kakko{\chi's}\cdot I_{i+k}^{-}
=\kakko{m\chi's}\cdot I_{i+k}^{-}=\kakko{\chi s}\cdot I_{i+k}^{-}$.
By induction we obtain the desired result. $\qed$

\subsubsection{calculation of the determinant of $S^\ast$}\label{seca.0.1}

By lemma \ref{lemmaa.1}, it follows that
$
E_1=x-a_{1,0}=x-\sum_{\chi\in\mathcal{X}_{M+K,0}}{\kakko{\chi}}
$,
$E_{k+1}=-a_{1,k}=-\sum_{\chi\in\mathcal{X}_{M+K,k}}{\kakko{\chi}}$.
Then,
\begin{align*}
\textstyle
\det{S^\ast}&\textstyle
=(-1)^{M+K+1}\cdot E_1=(-1)^{M+K}(\sum_{\chi\in\mathcal{X}_{M+K,0}}{\kakko{\chi}}-x)\\
&=(-1)^{M+K}(V_1^{t-(K-1)M}\cdots V_1^{t-M}V_1^tI_1^tI_1^{t-K}\cdots I_1^{t-(M-1)K}-x)\\
&= (-1)^{M+K}(U_1-x).
\end{align*}
$\qed$

\subsubsection{calculation of the determinant of $R^\ast$}\label{seca.0.2}

Let $w_k:=\sum_{\chi\in\mathcal{X}_{M+K-1,k}}{\kakko{\chi s}}$ and $w_{-1}:=0$.
By lemma \ref{lemmaa.2}, it follows that
\begin{align*}
\textstyle
a_{i,k}
&=\textstyle
\sum_{\chi\in\mathcal{X}_{M+K-1,k-1}}{\kakko{\chi s}}+
\sum_{\chi\in\mathcal{X}_{M+K-1,k}}{\kakko{\chi m}}\\
&=\textstyle
\sum_{\chi\in\mathcal{X}_{M+K-1,k-1}}{\kakko{\chi s}}+
\sum_{\chi\in\mathcal{X}_{M+K-1,k}}{\kakko{\chi s}}\cdot I_{k+1}^-\ \ 
=w_{k-1}+w_k\cdot I_{k+1}^-.
\end{align*}
Let $z_n:=I_1^-I_2^-\cdots I_n^-$.
By the cofactor expansion w.r.t.{}\,the $(M+K)$-th row, we obtain:
\begin{align*}
\det{R^\ast}&=(-1)^{M+K+1}
\{E_1-E_2z_1+E_3z_2-\cdots\\
&\hspace{3cm}
+(-1)^{M+K+1}E_{M+K}z_{M+K-1}
+(-1)^{M+K+1}z_{M+K}\}\\
&=(-1)^{M+K+1}
\{x-a_{1,0}+a_{1,1}z_1-a_{1,2}z_2+\cdots\\
&\hspace{2.5cm}+(-1)^{M+K}a_{1,M+K-1}\cdot z_{M+K-1}+(-1)^{M+K+1}z_{M+K}
\}\\
&=(-1)^{M+K+1}
\{
x-w_0I_1^-+w_0z_1-w_1z_1I_2^-+w_1z_2-w_2z_2I_3^-+\cdots\\
&\hspace{2cm}
+(-1)^{M+K}w_{M+K-1}z_{M+K-1}I_{M+K}^-+(-1)^{M+K+1}z_{M+K}\}\\
&=(-1)^{M+K+1}\cdot x.
\end{align*}
$\qed$

\subsubsection{calculation of the determinant of $L^\ast$}

Recall
the alternative form of the matrix $X_t(y)$ (cf.{}\,(\ref{eq2.2a})):
\begin{align*}
X_t(y)&=
R_{t-MK}(y)R_{t-(M+1)K}(y)\cdots R_{t-(2M-1)K}(y)\times\nonumber\\
&\hspace{2cm}
L_{t-(2N-1)M}(y)\cdots L_{t-(K+1)M}(y)L_{t-KM}(y),
\end{align*}
in which the rightmost matrix $L_{t-KM}(y)$ is essential.
The formula $\det{L^\ast}\!=\!(-1)^{M+K+1}\cdot x$ 
comes about through arguments similar to those in
section \ref{seca.0.2} concerning (\ref{eq2.2a}).
$\qed$


\begin{thebibliography}{10}
\bibitem{iwao} Iwao S 2008
``Solution of the generalised periodic Toda equation"
\textit{J. Phys. A. Math. Theor.} {\bf 41} 115201
\bibitem{iwao2} Iwao S 
``Solution of the generalised periodic Toda equation II;
Tau function solution" in submition to \textit{J.Phys.A.Math.Theor.},
(arXiv/0912.2213)
\bibitem{Tokihiro} Mada J, Idzumi M and Tokihiro T 2004
``Conserved quantities of generalized periodic box-ball
systems constructed from the ndKP equation"
\textit{J. Phys. A: Math. Gen.} {\bf 37} 6531--6556
\bibitem{Mumford} van Moerbeke P and Mumford D 1979
``The spectrum of difference operators and algebraic curves"
\textit{Acta Math.} {\bf 143} (1--2) 94--154
\bibitem{tata} Mumford D, Musili C, Nori M, Previato E and Stillman M
1983
\textit{Tata Lectures on Theta I} 
(\textit{Progress in mathematics}; v.28)
ed.{}\, 
Bass H, Oesterl\'{e} J and Weinstein A
(Berlin: Birkh\"{a}user)
\bibitem{Willox}
Willox R, Tokihiro T and Satsuma J 1997 
``Darboux and binary Darboux transformations for the nonautonomous
discrete KP equation" 
\textit{J. Math. Phys.} {\bf 38} 6455--6469
\end{thebibliography}
\end{document}